\documentclass[11pt,twoside,reqno,letter]{amsart}
\usepackage[margin=1in]{geometry}
\usepackage{amsmath,graphicx,varioref,amscd,amssymb,color,bm,amsthm,epsfig,color,enumerate,fancybox,bbm,esint,upref,xcolor,latexsym,mathtools,breqn,appendix,setspace,physics,tikz}
\usepackage[pdftex, colorlinks=true, linktocpage=true]{hyperref}

\DeclareUnicodeCharacter{221E}{$\infty$}

\makeatletter
\def\@tocline#1#2#3#4#5#6#7{\relax
  \ifnum #1>\c@tocdepth 
  \else
    \par \addpenalty\@secpenalty\addvspace{#2}%
    \begingroup \hyphenpenalty\@M
    \@ifempty{#4}{%
      \@tempdima\csname r@tocindent\number#1\endcsname\relax
    }{%
      \@tempdima#4\relax
    }%
    \parindent\z@ \leftskip#3\relax \advance\leftskip\@tempdima\relax
    \rightskip\@pnumwidth plus4em \parfillskip-\@pnumwidth
    #5\leavevmode\hskip-\@tempdima
      \ifcase #1
      \or\or \hskip 1em \or \hskip 2em \else \hskip 3em \fi%
      #6\nobreak\relax
    \dotfill\hbox to\@pnumwidth{\@tocpagenum{#7}}\par
    \nobreak
    \endgroup
  \fi}
  
  \def\l@subsection{\@tocline{2}{0pt}{2pc}{5pc}{}}
\makeatother

\numberwithin{equation}{section}
\newtheorem{thm}{Theorem}[section]
\newtheorem*{thm*}{Theorem}
\newtheorem{mydef}[thm]{Definition}
\newtheorem{lem}[thm]{Lemma}
\newtheorem*{lem*}{Lemma}

\newtheorem*{prop*}{Proposition}
\newtheorem{conj}[thm]{Conjecture}
\newtheorem{rem}[thm]{Remark}
\newcommand{\R}{\mathbb{R}}

\newcommand*\circled[1]{\tikz[baseline=(char.base)]{\node[shape=rectangle,draw,inner sep=2pt] (char) {#1};}}

\begin{document}


\title[Uniqueness in an NSE-NLS model of superfluids]{Uniqueness in a Navier-Stokes-nonlinear-Schr\"odinger model of superfluidity}

\author[Jayanti]{Pranava Chaitanya Jayanti}
\address[Jayanti]{\newline
Department of Physics \\ University of Maryland \\ College Park, MD 20742, USA.}
\email[]{\href{jayantip@umd.edu}{jayantip@umd.edu}}

\author[Trivisa]{Konstantina Trivisa}
\address[Trivisa]{\newline
Department of Mathematics \\ University of Maryland \\ College Park, MD 20742, USA.}
\email[]{\href{trivisa@math.umd.edu}{trivisa@math.umd.edu}}

\date{\today}



\thanks{P.C.J. was partially supported by the Ann Wylie Fellowship at UMD. Both P.C.J. and K.T. gratefully acknowledge the support of the National Science Foundation under the awards DMS-1614964 and DMS-2008568.}

\maketitle

\begin{abstract}
In \cite{Jayanti2022LocalSuperfluidity}, the authors proved the existence of local-in-time weak solutions to a model of superfluidity. The system of governing equations was derived in \cite{Pitaevskii1959PhenomenologicalPoint} and couples the nonlinear Schr\"odinger equation (NLS) and the Navier-Stokes equations (NSE). In this article, we prove a weak-strong type uniqueness theorem for these weak solutions. Only some of their regularity properties are used, allowing room for improved existence theorems in the future, with compatible uniqueness results. 
\end{abstract}

\setcounter{tocdepth}{2} 
\tableofcontents


\section{Introduction}  \label{intro}

This article deals with the problem of uniqueness of the local-in-time weak solutions to a model of superfluidity governed by the nonlinear Schr\"odinger equation (NLS) coupled with the inhomogeneous Navier-Stokes equations (NSE) for incompressible fluids. The NLS is used to describe the dynamics of the superfluid phase, while the NSE is employed to describe the evolution of the normal Helium liquid. \medskip

We present a result that is analogous to weak-strong uniqueness, but for less regular strong solutions (since our weak solutions are more regular than usual). Motivated by the analytical results in \cite{Jayanti2022LocalSuperfluidity}, we investigate the issue of uniqueness in Section \ref{local uniqueness of weak solutions} presenting a new class of weak solutions with additional regularity properties. Weak solutions usually lack the regularity that is needed to make them unique, a notable exception being the incompressible NSE in 2D \cite{Lions1959Un2}. Two different solutions to this impediment have been established in the literature. One is the conditional regularity theorems, where the weak solutions are shown to be smooth as soon as they belong to a critical regularity class. There are many results in this direction, mostly devoted to the incompressible NSE $-$ see Caffarelli et al. \cite{Caffarelli1982PartialEquations}, Escauriaza et al. \cite{EscauriazaL2003L3-solutionsUniqueness}, Prodi \cite{Prodi1959UnNavier-Stokes}, Serrin \cite{Serrin1963TheEquations}, or more recently, Neustupa et al. \cite{Neustupa2002AnVelocity,Neustupa2000AnEquations}. In the case of the compressible NSE, an analogous result was established in \cite{Feireisl2011SuitableFluids}. \medskip

The other approach, used in this work, is to show weak-strong uniqueness, i.e., that a strong solution is unique in the class of weak solutions (see \cite{Wiedemann2017Weak-strongDynamics} for a survey, including the case of the incompressible NSE). In the context of compressible fluids, related results can be found in Feireisl, Novotn\'{y} and collaborators \cite{Feireisl2011SuitableFluids,Feireisl2012RelativeSystem,Feireisl2012WeakstrongSystem}, and in Mellet and Vasseur \cite{Mellet2008AsymptoticEquations}. Intermediate on the spectrum of complexity associated with compressibility are the incompressible inhomogeneous fluids, which is the model considered in this paper. There has been work done for weak-strong uniqueness in this regard too \cite{Germain2008StrongSystem}. Similar methods have been used to obtain conditional uniqueness for related systems, like fluid-solid interactions \cite{Chemetov2019Weak-strongCondition,Kreml2020Weak-strongInteraction}, fluid-particle interactions \cite{Ballew2013WeaklySystem}, compressible magnetohydrodynamics \cite{Yan2013OnFlows}, and Euler-Korteweg-Poisson systems \cite{Donatelli2015Well/IllProblems}. \medskip


Our analysis is motivated by the pioneering work of Dafermos \cite{Dafermos1979TheStability} and DiPerna \cite{DiPerna1979UniquenessLaws}, the results of Germain \cite{Germain2011WeakstrongSystem}, the analysis of Mellet and Vasseur \cite{Mellet2008AsymptoticEquations}, as well as the approach of Feireisl et al. \cite{Feireisl2011SuitableFluids}. Usually, weak-strong uniqueness theorems involve strong solutions that are very regular: bounded (in $x$) first-order space-time derivatives. This is required so that the strong solution bears the brunt of the smoothness requirements, allowing us to work with rather irregular weak solutions (for instance, those satisfying a simple $L^2_x$ based energy estimate). A Gronwall’s argument is employed on the difference between the weak and strong solutions to conclude uniqueness (see Theorem \ref{weak-strong uniqueness} below), yielding that a weak solution agrees with a classical solution starting from the same initial data when such a classical solution exists. We also refer the reader to Section 2.5 in \cite{Lions1996MathematicalMechanics} where a related result is presented for incompressible inhomogeneous fluids. It must be noted that the scheme, originally developed for weak solutions of the Navier-Stokes equations (with regularity only at the level of the energy bound, i.e., $u\in L^2_x$), does not work here because of added interactions with the Schr\"odinger equation. We will necessarily need to consider solutions with smoother velocity fields. We will discuss this briefly in Section \ref{weak-strong uniqueness proof} below, after sketching the steps used in the treatment of the momentum equation. \medskip

The weak solutions constructed in \cite{Jayanti2022LocalSuperfluidity} to a model of superfluidity required more a priori estimates than just the energy, and this led to better-behaved solutions (albeit not regular enough to prove uniqueness). We can, however, trade the increased regularity of our weak solution for a less smooth strong solution. A \textit{weak-moderate} uniqueness theorem is presented in Section \ref{local uniqueness of weak solutions}. Apart from the use of ``moderate" solutions (which belong to a lower regularity class), another feature is that the full regularity of the existence result is not utilized. This means that this uniqueness result will continue to hold even if weaker solutions (than in \cite{Jayanti2022LocalSuperfluidity}) are shown to exist for the model in question. \medskip

\subsection{Notation} \label{notation}
Let $\mathfrak{D}(\Omega)$ be the space of smooth, compactly-supported functions on $\Omega$. Then, $H^s_0(\Omega)$ is the completion of $\mathfrak{D}$ under the Sobolev norm $H^s$. The more general Sobolev spaces are denoted by $W^{s,p}(\Omega)$, where $s\in\R$ is the derivative index and $1\le p\le\infty$ is the integrability index. A dot on top, like $\Dot{H}^s(\Omega)$ or $\Dot{W}^{s,p}$, is used when referring to the homogeneous Sobolev spaces. \medskip

Consider a 3D vector-valued function $u\equiv (u_1,u_2,u_3)$, where $u_i\in \mathfrak{D}(\Omega), i=1,2,3$. The collection of all divergence-free functions $u$ defines $\mathfrak{D}_d(\Omega)$. Then, $H^s_d(\Omega)$ is the completion of $\mathfrak{D}_d(\Omega)$ under the $H^s$ norm. In addition, to say that a complex-valued wavefunction $\psi\in H^s(\Omega)$ means that its real and imaginary parts are the limits (in the $H^s$ norm) of functions in $\mathfrak{D}(\Omega)$. \medskip





The $L^2$ inner product, denoted by $\langle \cdot,\cdot \rangle$, is sesquilinear (the first argument is complex conjugated, indicated by an overbar) to accommodate the complex nature of the Schr\"odinger equation. Thus, for example, $\langle \psi,B\psi \rangle = \int_{\Omega} \Bar{\psi}B\psi \ dx$. Needless to say, since the velocity and density are real-valued functions, we will ignore the complex conjugation when they constitute the first argument of the inner product. \medskip

We use the subscript $x$ on a Banach space to denote the Banach space is defined over $\Omega$. For instance, $L^p_x$ stands for the Lebesgue space $L^p(\Omega)$, and similarly for the Sobolev spaces: $H^s_{d,x} := H^s_d(\Omega)$. For spaces/norms over time, the subscript $t$ will denote the time interval in consideration, such as $L^p_t := L^p_{[0,T]}$, where $T$ stands for the local existence time unless mentioned otherwise. The Bochner spaces $L^p(0,T;X)$ and $C([0,T];X)$ have their usual meanings, as ($L^p$ and continuous, respectively) maps from $[0,T]$ to a Banach space $X$. 
\medskip

We also use the notation $X\lesssim Y$ to imply that there exists a positive constant $C$ such that $X\le CY$. The dependence of the constant on various parameters (including the initial data), will be denoted using a subscript as $X\lesssim_{k_1,k_2} Y$ or $X\le C_{k_1,k_2}Y$. \bigskip

\subsection{Organization of the paper}

In Section \ref{mathematical model}, we describe the mathematical model and recall the existence results of \cite{Jayanti2022LocalSuperfluidity}, followed by the statement of the main theorem. Section \ref{weak-strong uniqueness proof} is devoted to an explanation of how the standard weak-strong uniqueness is difficult to establish for this model unless the weak solutions are sufficiently regular (at which point weak-moderate uniqueness itself applies). Finally, the main result concerning the moderate solutions is proved in Section \ref{local uniqueness of weak solutions}, using the energy structure of the system. \bigskip

\section{Mathematical model and main results} \label{mathematical model}

\subsection{Model and recap of existence results}

The nonlinear Schr\"odinger equation (NLS) is used to describe the dynamics of the superfluid phase, while the incompressible inhomogeneous Navier-Stokes equations (NSE) govern the normal Helium liquid (see \cite{Jayanti2021GlobalEquations} for a discussion on the physics of superfluidity and other related models). The ``Pitaevskii model" that we consider in this work is as follows. Henceforth, all statements are in the context of a smooth, bounded domain $\Omega\subset\R^3$. \medskip

    
    
    

\begin{align}
    \partial_t \psi + \Lambda B\psi &= -\frac{1}{2i}\Delta\psi + \frac{\mu}{i}\lvert\psi\rvert^2 \psi \tag{NLS} \label{NLS} \\
    B = \frac{1}{2}\left(-i\nabla - u \right)^2 + \mu \lvert \psi \rvert^2 &= -\frac{1}{2}\Delta + iu\cdot\nabla + \frac{1}{2}\lvert u \rvert^2 + \mu \lvert \psi \rvert^2 \tag{CPL} \label{coupling} \\
    \partial_t \rho + \nabla\cdot(\rho u) &= 2\Lambda\Re(\Bar{\psi}B\psi) \tag{CON} \label{continuity} \\
    \partial_t (\rho u) + \nabla\cdot (\rho u \otimes u) + \nabla p - \nu \Delta u &= \!\begin{multlined}[t]
    -2\Lambda \Im(\nabla \Bar{\psi}B\psi) + \Lambda\nabla \Im(\Bar{\psi}B\psi) \tag{NSE} \label{NSE} + \frac{\mu}{2}\nabla\lvert\psi\rvert^4 
    \end{multlined} \\
    \nabla\cdot u &= 0 \tag{DIV} \label{divergence-free}
\end{align}\medskip

We supplement the equations with initial and boundary conditions\footnote{For a justification of the exclusion of $t=0$ in the boundary conditions for the wavefunction, see Remark 2.4 in \cite{Jayanti2022LocalSuperfluidity}.} as follows.

\begin{equation*} \tag{INI} \label{initial conditions}
        \psi(0,x) = \psi_0(x) \qquad u(0,x) = u_0(x) \qquad \rho(0,x) = \rho_0(x) \quad a.e. \ x\in\Omega
\end{equation*}\medskip

\begin{equation} \tag{BC} \label{boundary conditions}
    \begin{aligned}
    u = \frac{\partial u}{\partial n} = 0 \quad &a.e. \ (t,x)\in[0,T]\times\partial\Omega \\
    \psi = \frac{\partial \psi}{\partial n} = \frac{\partial^2 \psi}{\partial n^2} = \frac{\partial^3 \psi}{\partial n^3} = 0 \quad &a.e. \ (t,x)\in(0,T]\times\partial\Omega \notag
    \end{aligned}
\end{equation}\medskip

\noindent where $n$ is the outward normal direction on the boundary $\partial\Omega$, and $T$ is the local existence time. For more details about the model, refer to the discussion in \cite{Jayanti2022LocalSuperfluidity} or the derivation in \cite{Pitaevskii1959PhenomenologicalPoint}. Weak solutions that we seek are those that satisfy the governing equations in the sense of distributions, for a certain class of test functions. \medskip

\begin{mydef}[Weak solutions\footnote{See Remark \ref{strong or weak solutions?}.}] \label{definition of weak solutions}
    Let $\Omega \subset \R^3$ be a bounded set with a smooth boundary $\partial\Omega$. For a given time $T>0$, consider the following test functions:
    \begin{enumerate}
        \item a complex-valued scalar field $\varphi \in H^1(0,T;L^2(\Omega))\cap L^2(0,T;H_0^1(\Omega))$,
        
        \item a real-valued, divergence-free (3D) vector field $\Phi \in H^1(0,T;L^2_d(\Omega)) \cap L^2(0,T;H^1_d(\Omega))$, and
        
        \item a real-valued scalar field $\sigma \in H^1(0,T;L^2(\Omega))\cap L^2(0,T;H^1(\Omega))$.
    \end{enumerate}
    
    A triplet $(\psi,u,\rho)$ is called a \textbf{weak solution} to the Pitaevskii model if:
    
    \begin{enumerate} [(i)]
        \item 
        
        \begin{equation}
        \begin{gathered}
            \psi\in L^2(0,T;H^{\frac{7}{2}+\delta}_0(\Omega)) \\
            u\in L^2(0,T;H^{\frac{3}{2}+\delta}_d(\Omega)) \\
            \rho \in L^{\infty}([0,T]\times\Omega)
        \end{gathered}
        \end{equation}\medskip
        
        \item and they satisfy the governing equations in the sense of distributions for all test functions, i.e.,
        
        \begin{multline} \label{weak solution wavefunction}
            -\int_0^T \int_{\Omega} \left[ \psi\partial_t\Bar{\varphi} + \frac{1}{2i}\nabla\psi\cdot\nabla\Bar{\varphi} - \Lambda\Bar{\varphi}B\psi - i\mu\Bar{\varphi}\lvert \psi \rvert^2\psi \right] dx \ dt \\ = \int_{\Omega} \left[ \psi_0\Bar{\varphi}(t=0) - \psi(T)\Bar{\varphi}(T) \right] dx
        \end{multline}
    
        \begin{multline} \label{weak solution velocity}
            -\int_0^T \int_{\Omega} \left[ \rho u\cdot \partial_t \Phi + \rho u\otimes u:\nabla\Phi - \nu\nabla u:\nabla\Phi - 2\Lambda\Phi\cdot\Im(\nabla\Bar{\psi}B\psi) \right] dx \ dt \\ = \int_{\Omega} \left[ \rho_0 u_0 \Phi(t=0) - \rho(T)u(T)\Phi(T) \right] dx
        \end{multline}
    
        \begin{equation} \label{weak solution density}
            -\int_0^T \int_{\Omega} \left[ \rho \partial_t \sigma + \rho u\cdot\nabla\sigma + 2\Lambda \sigma \Re(\Bar{\psi}B\psi) \right] dx \ dt = \int_{\Omega} \left[ \rho_0 \sigma(t=0) - \rho(T)\sigma(T) \right] dx
        \end{equation}
    
        \noindent where (the initial data) $\psi_0 \in H^{\frac{5}{2}+\delta}_0(\Omega)$, $u_0 \in H^{\frac{3}{2}+\delta}_d(\Omega)$ and $\rho_0 \in L^{\infty}(\Omega)$. \medskip
    \end{enumerate}\medskip
    
\end{mydef}\medskip



    

In \cite{Jayanti2022LocalSuperfluidity}, the following existence result was proven. \medskip

\begin{thm} [Local existence] \label{local existence}
    For any $\delta \in (0,\frac{1}{2})$, let $\psi_0 \in H^{\frac{5}{2}+\delta}_0(\Omega)$ and $u_0\in H^{\frac{3}{2}+\delta}_d(\Omega)$. Suppose $\rho_0$ is bounded both above and below a.e. in $\Omega$, i.e., $0<m\le\rho_0\le M<\infty$. Then, there exists a local existence time $T$ and at least one weak solution $(\psi,u,\rho)$ to the Pitaevskii model. In particular, the weak solutions have the following regularity:
    
    \begin{gather}
        \psi\in C([0,T];H^{\frac{5}{2}+\delta}_0(\Omega)) \cap L^2(0,T;H^{\frac{7}{2}+\delta}_0(\Omega)) \label{weak solution psi regularity} \\
        u\in C([0,T];H^{\frac{3}{2}+\delta}_d(\Omega)) \cap L^2(0,T;H^{2}_d(\Omega)) \label{weak solution u regularity} \\
        \rho\in L^{\infty}([0,T]\times \Omega)\cap C([0,T];L^2(\Omega)) \label{weak solution rho regularity}
    \end{gather}
    
    \noindent where $T$ depends (inversely) on:
    \begin{enumerate}[(i)]
        \item $\varepsilon\in (0,m)$, the allowed infimum of the density field, and
        \item $X_0 = 1 + \norm{\Delta\psi_0}_{L^2_x}^2 + \nu\norm{\nabla u_0}_{L^2_x}^2$
    \end{enumerate}
    
    
\end{thm}\medskip

\begin{rem} \label{strong or weak solutions?}
    The regularity of the solutions seem to suggest that the wavefunction and velocity are strong solutions. Indeed this is true, as they are strongly continuous in their topologies. On the other hand, the density is truly a weak solution and is the reason for referring to the triplet as a weak solution.
\end{rem} \medskip

\begin{rem}
    For fixed values of $\nu$ and $\Lambda$, it was argued in Section 2.1 of \cite{Jayanti2022LocalSuperfluidity} that the local existence time $T$ grows with decreasing $\varepsilon$ (up to a finite maximum value).
\end{rem}\medskip

The solutions also satisfy some bounds/estimates (some derived a priori), summarized here for convenience. \medskip 

\begin{enumerate}
    \item \textit{Superfluid mass bound}:
    \begin{equation} \label{superfluid mass bound}
    \norm{\psi}_{L^2_x}(t) \le \norm{\psi_0}_{L^2_x} \quad a.e. \  t \in [0,T]
    \end{equation}\medskip
    
    \item \textit{Normal fluid density bound}: Given our choice of $T$, we know that
    
    \begin{equation} \label{density upper lower bound}
        0<\varepsilon<\rho(t,x)<M' := M+m-\varepsilon \ , \ a.e. \ (t,x)\in [0,T]\times\Omega
    \end{equation}\medskip 
    
    \item \textit{Energy equation}:
    \begin{equation} \label{energy equality}
    \begin{multlined}
        \left( \frac{1}{2}\norm{\sqrt{\rho}u}_{L^2_x}^2 + \frac{1}{2}\norm{\nabla\psi}_{L^2_x}^2 + \frac{\mu}{2}\norm{\psi}_{L^4_x}^4 \right)(t) + \nu\norm{\nabla u}_{L^2_{[0,t]}L^2_x}^2 + 2\Lambda\norm{B\psi}_{L^2_{[0,t]}L^2_x}^2 \\ = \frac{1}{2}\norm{\sqrt{\rho_0}u_0}_{L^2_x}^2 + \frac{1}{2}\norm{\nabla\psi_0}_{L^2_x}^2 + \frac{\mu}{2}\norm{\psi_0}_{L^4_x}^4 =: E_0 \quad a.e. \ t\in [0,T]
    \end{multlined}
    \end{equation}\medskip
    
    \item \textit{Higher-order energy estimate}:
    \begin{equation} \label{higher-order energy estimate}
    \begin{gathered}
        X(t) \le 2X_0 \quad a.e. \ t\in[0,T] \qquad , \qquad \int_0^T Y(\tau) d\tau \le 31X_0 \\
        where \ X(t) = 1 + \norm{\Delta\psi}_{L^2_x}^2 + \nu\norm{\nabla u}_{L^2_x}^2 \quad , \quad X_0 = X(0) \\
        Y = \Lambda\norm{\nabla(B\psi)}_{L^2_x}^2 + \norm{\sqrt{\rho}\partial_t u}_{L^2_x}^2 + \frac{\nu^2}{M'}\norm{\Delta u}_{L^2_x}^2
    \end{gathered}
    \end{equation}\medskip
    
    \item \textit{Highest-order energy estimate}:
    \begin{equation} \label{highest-order energy estimate}
    \begin{gathered}
        \norm{\psi}_{H^{\frac{5}{2}+\delta}_x}^2 (t) \lesssim \norm{\psi_0}_{H^{\frac{5}{2}+\delta}_x}^2 e^{c Q_T} \quad a.e. \ t\in[0,T] \\
        \norm{B\psi}_{L^2_{[0,T]}H^{\frac{3}{2}+\delta}_x} \lesssim \Lambda^{-\frac{1}{2}} (Q_T^{\frac{1}{2}} + 1) e^{c Q_T} \norm{\psi_0}_{H^{\frac{5}{2}+\delta}_x} \\
        \norm{u}_{H^{\frac{3}{2}+\delta}_x}^2(t) \lesssim \norm{u_0}_{H^{\frac{3}{2}+\delta}_x}^2 + \frac{1}{\nu}\sqrt{\frac{M'}{\varepsilon}}X_0 \\
        where \ Q_T = \Lambda\frac{M'}{\nu^2}X_0 + \left(\Lambda \frac{M'}{\nu^2\varepsilon} + \gamma\right)X_0^2 T + \Lambda E_1^2 T \quad , \quad E_1 = \norm{u_0}_{H^{\frac{3}{2}+\delta}_x}^2 + \norm{\psi_0}_{H^{\frac{5}{2}+\delta}_x}^2
    \end{gathered}
    \end{equation}\medskip
    
    \item \textit{Time-derivative bounds}:
    \begin{equation} \label{time-derivative bounds}
    \begin{gathered}
        \norm{\partial_t \psi}_{L^2_{[0,T]}L^2_x} \lesssim \norm{B\psi}_{L^2_{[0,T]}L^2_x} + T^{\frac{1}{2}}\norm{\psi}_{L^{\infty}_{[0,T]}H^2_x} + \mu T^{\frac{1}{2}} \norm{\psi}_{L^{\infty}_{[0,T]}H^1_x}^3 \\
        \norm{\partial_t u}_{L^2_{[0,T]}L^2_x} \lesssim \frac{1}{\sqrt{\varepsilon}} X_0 \\
        \norm{\partial_t \rho}_{L^2_{[0,T]}H^{-1}_x} \lesssim \left( T E_0 M' \right)^{\frac{1}{2}} + \left( \Lambda^{-1} E_0 X_0 \right)^{\frac{1}{2}}
    \end{gathered}
    \end{equation}\medskip
    
\end{enumerate} \medskip

In the derivation of \eqref{energy equality}, we make use of the \textit{non-conservative form} of \eqref{NSE}. Note that $\Tilde{p}$ is simply a modified pressure that includes gradient terms on the RHS of the momentum equation. \medskip

\begin{equation}  \label{NSE'}
    \rho\partial_t u + \rho u\cdot\nabla u + \nabla \Tilde{p} - \nu \Delta u = -2\Lambda \Im(\nabla \Bar{\psi}B\psi) - 2\Lambda u\Re(\Bar{\psi}B\psi) \tag{NSE'}
\end{equation}\bigskip

\subsection{Uniqueness results}

Given that weak solutions exist, it is instructive to ask whether they are unique and under what additional assumptions, if any. The main result of this article is a ``weak-moderate" uniqueness theorem. But first, we will discuss the more standard (and desirable) weak-strong uniqueness statement, i.e., starting from the same initial conditions, if we have a weak solution and a strong solution (in the precise sense described below), then they are both identical almost everywhere. \medskip

\begin{conj} [Weak-strong uniqueness] \label{weak-strong uniqueness}
    Let $(\psi,u,\rho)$ be a weak solution to the Pitaevskii model, whose regularity is governed by \eqref{superfluid mass bound}-\eqref{energy equality}, i.e., 
    
    \begin{gather}
        \psi\in C([0,T];H^1_0(\Omega)) \cap L^2(0,T;H^2_0(\Omega)) \label{weak solution psi regularity least} \\
        u\in L^{\infty}([0,T];L^2_d(\Omega)) \cap L^2(0,T;H^1_d(\Omega)) \label{weak solution u regularity least} \\
        \rho\in L^{\infty}([0,T]\times \Omega)\cap C([0,T];L^2(\Omega)) \label{weak solution rho regularity least}
    \end{gather}
    
    Suppose there exists a strong solution $(\Tilde{\psi},\Tilde{u},\Tilde{\rho})$ starting from the same initial conditions in $H^1_0 \times L^2_d \times L^{\infty}$ (such that both these solutions exist up to some time $T$). Specifically, it is sufficient that the strong solution has the following additional regularity in addition to \eqref{weak solution psi regularity least}-\eqref{weak solution rho regularity least}:
    
    \begin{equation} \label{increased regularity for weak-strong uniqueness}
    \begin{gathered}
        \Tilde{\rho}\in L^2(0,T;W^{1,3}(\Omega)) \quad , \quad \Tilde{u} \in L^2(0,T;W^{1,\infty}(\Omega)) \\
        \partial_t \Tilde{u} \in L^2(0,T;L^{\infty}(\Omega))
    \end{gathered}
    \end{equation}\medskip
    
    Then, the solutions are identical a.e. in $[0,T]\times\Omega$.
\end{conj}\medskip

\begin{rem} \label{global weak solutions}
    If we have weak solutions whose regularity is determined purely by the energy equation, then it is tempting to assume that such solutions would be global in time. Unfortunately, this is not the case, since we do not have a way of guaranteeing a positive lower bound for the density, globally in time. At the moment, it is not clear at all how to approach such a problem, due to the presence of the sign-indefinite term on the RHS of \eqref{continuity}. A more elaborate discussion on this issue can be found in \cite{Jayanti2022LocalSuperfluidity}, and this is an interesting open question worth exploring in the future.
\end{rem}\medskip

Based on the regularity in \eqref{weak solution psi regularity least}, one can easily see that $\partial_t \psi$ is already in $L^2(0,T;L^2(\Omega))$, which is why no improvement is needed for the smoothness of $\Tilde{\psi}$. In fact, the wavefunction has enough regularity to be classified as a strong solution; it is the velocity and density that make the solutions ``weak". While this conjecture works for the Navier-Stokes equations, we run into a problem here: the coupling to the NLS. This gives rise to terms that cannot be handled at the level of regularity of the weak solutions. To this end, we will outline the main steps in the analysis of the momentum equation, to clarify exactly where the difficulties arise. \medskip

Furthermore, in the above conjecture, we observe that $\Tilde{u}$ needs to have bounded (in space) spatio-temporal derivatives. As much as this is not unusual to expect from strong solutions, it is possible to trade this rather high regularity for slightly smoother weak solutions. This is stated in the next theorem, and will be referred to as ``weak-moderate uniqueness" since the notion of strong solution is significantly relaxed. \medskip

\begin{thm} [Weak-moderate uniqueness] \label{local uniqueness}
    Let $(\psi,u,\rho)$ be a weak solution to the Pitaevskii model, with regularity governed by \eqref{superfluid mass bound}-\eqref{higher-order energy estimate}, i.e.,
    
    \begin{gather}
        \psi\in C([0,T];H^{2}_0(\Omega)) \cap L^2(0,T;H^{3}_0(\Omega)) \label{weak solution psi regularity reduced} \\
        u\in C([0,T];H^{1}_d(\Omega)) \cap L^2(0,T;H^{2}_d(\Omega)) \label{weak solution u regularity reduced} \\
        \rho\in L^{\infty}([0,T]\times \Omega)\cap C([0,T];L^2(\Omega)) \label{weak solution rho regularity reduced}
    \end{gather}
    
    Suppose there exists a ``moderate" solution $(\Tilde{\psi},\Tilde{u},\Tilde{\rho})$ starting from the same initial data (in $H^2_0 \times H^1_d \times L^{\infty}$), such that $\Tilde{\rho} \in L^2(0,T;W^{1,3}(\Omega)), \ \partial_t\Tilde{u}\in L^2(0,T;L^3(\Omega))$ in addition to the regularity in \eqref{weak solution psi regularity reduced}-\eqref{weak solution rho regularity reduced}. Then, the solutions are identical a.e. in $[0,T]\times\Omega$. \medskip 
    
        
        
        
        
        
 \end{thm}\medskip


It is to be noted that in the proof of Theorem \ref{local existence}, we already established that for a weak solution, $\partial_t u\in L^2(0,T;L^2)$ $-$ evident from \eqref{density upper lower bound} and \eqref{higher-order energy estimate}. So, as far as the (time derivative of the) ``moderate" velocity goes, we have only assumed a slight increase in the spatial integrability, from $L^2$ to $L^3$ (as opposed to $L^{\infty}$ in weak-strong uniqueness). Also noteworthy is the absence of any additional (over the weak solution) bound on the spatial derivatives of $\Tilde{u}$. The regularity of the ``moderate" density remains unchanged from Conjecture \ref{weak-strong uniqueness}. \medskip

\begin{rem} \label{lower regularity of weak solutions in weak-moderate result}
    It is important to recognize that the regularity of the weak solutions in Theorem \ref{local uniqueness} is lower than those in Theorem \ref{local existence}. As mentioned before, this allows wiggle room for enhanced (less regular) existence results without compromising uniqueness. 
\end{rem}\bigskip

\bigskip

Summarized below are some well-known properties of Sobolev and Lebesgue spaces (see \cite{Temam1977Navier-StokesAnalysis} and Chapter 5 of \cite{Evans2010PartialEquations}) that will be repeatedly utilized in the calculations that follow. We will begin with a discussion on the difficulties of proving Conjecture \ref{weak-strong uniqueness} in the next section. \medskip

\begin{lem}[Sobolev embeddings] \label{sobolev embeddings}
    For $\Omega$ a smooth, bounded subset of $\R^3$,
    \begin{enumerate}
        \item $H^1(\Omega)\subset L^6(\Omega)$ 
        \item $H^s(\Omega) \subset L^{\infty}(\Omega) \ \forall \ s>\frac{3}{2}$ \ ; \ $H^s(\R) \subset L^{\infty}(\R) \ \forall \ s>\frac{1}{2}$
        \item $H^s(\Omega) \Subset H^{s'}(\Omega) \ \forall \ s,s'\in\R, s>s'$
    \end{enumerate}
\end{lem}\medskip

\begin{lem}[Inequalities and interpolations] \label{inequalities and interpolations}
    For $\Omega$ a smooth, bounded subset of $\R^3$,
    \begin{enumerate}
        \item Poincar\'e's inequality: $\norm{f}_{L^p} \lesssim \norm{\nabla f}_{L^p} \ \forall  \ 1<p<\infty, f\in W^{1,p}_0(\Omega) \Rightarrow \norm{f}_{W^{1,p}} \equiv \norm{\nabla f}_{L^p}$
        \item Ladyzhenskaya's inequality: $\norm{f}_{L^4} \lesssim \norm{f}_{L^2}^{\frac{1}{4}} \norm{\nabla f}_{L^2}^{\frac{3}{4}} \ \forall \ f\in H^1_0(\Omega)$
        \item Agmon's inequality: $\norm{f}_{L^{\infty}} \lesssim \norm{f}_{H^1}^{\frac{1}{2}} \norm{f}_{H^2}^{\frac{1}{2}} \ \forall \ f\in H^1_0(\Omega)\cap H^2(\Omega)$
        \item Lebesgue interpolation: $\norm{f}_{L^3} \lesssim \norm{f}_{L^2}^{\frac{1}{2}} \norm{f}_{L^6}^{\frac{1}{2}} \ \forall \ f\in L^2 \cap L^6$
    \end{enumerate}
\end{lem}\bigskip

\section{Remark on weak-strong uniqueness (Conjecture \ref{weak-strong uniqueness})} \label{weak-strong uniqueness proof}

In this section, we will investigate Conjecture \ref{weak-strong uniqueness}, and our modus operandi is motivated by the classical results of Section 2.5 in \cite{Lions1996MathematicalMechanics}. Our goal is to demonstrate that the regularity in \eqref{weak solution psi regularity least}-\eqref{weak solution rho regularity least} is not sufficient to show a weak-strong uniqueness result, even if the strong solutions are smooth. \medskip

Consider the momentum equation. In what follows, $\Psi = -2\Lambda u \Re(\overline{\psi}B\psi)$ and $\Psi' = -2\Lambda\Im(\nabla\overline{\psi}B\psi)$. We begin from \eqref{weak solution velocity}, setting $\Phi = \Tilde{u}$. This way, we can arrive at an expression for $\int_{x} \rho u\cdot\Tilde{u}$. \medskip

\begin{equation} \label{rho u.Tilde(u) 1}
    \int_{\Omega} \rho u\cdot\Tilde{u} + \nu\int_{t,\Omega} \nabla u:\nabla\Tilde{u} = \int_{\Omega} \rho_0 \abs{u_0}^2 + \int_{t,\Omega} \rho u\cdot\left( \partial_t \Tilde{u} + \Tilde{u}\cdot\nabla\Tilde{u} \right) + \int_{t,\Omega} \Psi' \cdot\Tilde{u}
\end{equation}\medskip

We write \eqref{NSE'} for the strong solution, and rearrange it to get:

\begin{equation} \label{NSE strong solution - rearranged}
    \rho\partial_t \Tilde{u} + \rho u\cdot\nabla\Tilde{u} + \nabla \Tilde{p} -\nu\Delta\Tilde{u} = \Tilde{\Psi} + \Tilde{\Psi'} + (\rho - \Tilde{\rho})\partial_t \Tilde{u} + (\rho - \Tilde{\rho})\Tilde{u}\cdot\nabla\Tilde{u} + \rho (u - \Tilde{u})\cdot\nabla\Tilde{u}
\end{equation}\medskip

Multiplying \eqref{NSE strong solution - rearranged} by $u$ and integrating, \medskip

\begin{equation} \label{NSE strong solution - rearranged, multiplied by u}
\begin{multlined}
    \int_{t,\Omega}\rho u\partial_t \Tilde{u} + \int_{t,\Omega}\rho u\otimes u:\nabla\Tilde{u} +\nu\int_{t,\Omega}\nabla u:\nabla\Tilde{u} = \int_{t,\Omega}\Tilde{\Psi}\cdot u + \int_{t,\Omega}\Tilde{\Psi'}\cdot u \\ + \int_{t,\Omega}(\rho - \Tilde{\rho})u\cdot\left(\partial_t \Tilde{u} + \Tilde{u}\cdot\nabla\Tilde{u}\right) + \int_{t,\Omega}\rho (u - \Tilde{u})\otimes u:\nabla\Tilde{u}
\end{multlined}
\end{equation}\medskip

Adding \eqref{rho u.Tilde(u) 1} and \eqref{NSE strong solution - rearranged, multiplied by u}, \medskip

\begin{equation} \label{rho u.Tilde(u)}
\begin{multlined}
    \int_{\Omega} \rho u\cdot\Tilde{u} + 2\nu\int_{t,\Omega} \nabla u:\nabla\Tilde{u} = \int_{\Omega} \rho_0 \abs{u_0}^2 + \int_{t,\Omega}(\rho - \Tilde{\rho})u\cdot\left(\partial_t \Tilde{u} + \Tilde{u}\cdot\nabla\Tilde{u}\right) \\ + \int_{t,\Omega}\rho (u - \Tilde{u})\otimes u:\nabla\Tilde{u} + \int_{t,\Omega} \left[ \Psi' \cdot\Tilde{u} + \Tilde{\Psi} \cdot u + \Tilde{\Psi'} \cdot u \right]
\end{multlined}
\end{equation}\medskip

Next, we take the inner product of both \eqref{NSE} and \eqref{NSE'} with $u$, use incompressibility, and add them, to arrive at the energy equation for the normal fluid alone.

\begin{equation} \label{normal fluid energy equation}
    \frac{d}{dt}\frac{1}{2}\norm{\sqrt{\rho}u}_{L^2_x}^2 + \nu\norm{\nabla u}_{L^2_x}^2 = - 2\Lambda\int_{\Omega}u\cdot\Im(\nabla \Bar{\psi}B\psi) - \Lambda\int_{\Omega}\abs{u}^2\Re(\Bar{\psi}B\psi)
\end{equation}\medskip

This obviously holds as an a priori estimate, valid for the approximate fields $(\psi^N,u^N,\rho^N)$, but given the regularity of the solution, we can easily pass to the limit to see that it is accurate for the weak solutions too. From this equation, integrating on $[0,t]$ for $0\le t\le T$, we can obtain an equation for $\int_{\Omega} \rho u\cdot u$. \medskip

\begin{equation} \label{rho u.u}
    \frac{1}{2}\int_{\Omega} \rho \abs{u}^2 + \nu\int_{t,\Omega} \abs{\nabla u}^2 = \frac{1}{2}\int_{\Omega} \rho_0 \abs{u_0}^2 + \int_{t,\Omega} \left[ \frac{1}{2}\Psi\cdot u + \Psi' \cdot u \right]
\end{equation}\medskip

Finally, multiply \eqref{NSE strong solution - rearranged} by $\Tilde{u}$ and integrate over $[0,t]\times \Omega$. Since $\Tilde{u}$ is a strong solution, $\frac{1}{2}\abs{\Tilde{u}}^2$ can work as a test function. Therefore, we use \eqref{weak solution density} to simplify, and the resulting equation is: \medskip

\begin{equation} \label{rho Tilde(u).Tilde(u)}
\begin{multlined}
    \frac{1}{2}\int_{\Omega} \rho \abs{\Tilde{u}}^2 + \nu\int_{t,\Omega} \abs{\nabla \Tilde{u}}^2 = \frac{1}{2}\int_{\Omega} \rho_0 \abs{u_0}^2 + \int_{t,\Omega}(\rho-\Tilde{\rho})\Tilde{u}\cdot\left( \partial_t \Tilde{u} + \Tilde{u}\cdot\nabla\Tilde{u} \right) \\ + \int_{t,\Omega}\rho (u-\Tilde{u})\otimes \Tilde{u}:\nabla\Tilde{u} + \int_{t,\Omega} \left[ \Tilde{\Psi} + \Tilde{\Psi'} \right]\cdot\Tilde{u} + \Lambda\Re\int_{t,\Omega}\overline{\psi}B\psi \abs{\Tilde{u}}^2
\end{multlined}
\end{equation}\medskip

Add \eqref{rho u.u} and \eqref{rho Tilde(u).Tilde(u)}, and subtract \eqref{rho u.Tilde(u)}. Then, differentiate the result with respect to $t$. \medskip


\begin{equation} \label{weak-strong uniqueness for velocity - preliminary}
    \begin{aligned}
        \renewcommand{\arraystretch}{2}
        \left. \begin{array}{c}
         \frac{d}{dt}\frac{1}{2}\norm{\sqrt{\rho}(u - \Tilde{u})}_{L^2_x}^2  \\
          + \ \nu\norm{\nabla(u-\Tilde{u})}_{L^2_x}^2 
        \end{array}  \right\} = &-\underbrace{\int_{\Omega} (\rho-\Tilde{\rho})(u-\Tilde{u})\cdot\left( \partial_t \Tilde{u} + \Tilde{u}\cdot\nabla\Tilde{u} \right)}_{\circled{I-1}} - \underbrace{\int_{\Omega} \rho(u-\Tilde{u})\otimes(u-\Tilde{u}):\nabla\Tilde{u}}_{\circled{I-2}} \\
        &- \underbrace{\Lambda\Re\int_{\Omega} (\overline{\psi}B\psi - \overline{\Tilde{\psi}}\Tilde{B}\Tilde{\psi})(\abs{u}^2 - \abs{\Tilde{u}}^2)}_{\circled{I-3}} - \underbrace{\Lambda\Re\int_{\Omega} \overline{\Tilde{\psi}}\Tilde{B}\Tilde{\psi}\abs{u-\Tilde{u}}^2}_{\circled{I-4}} \\
        &- \underbrace{2\Lambda\Im\int_{\Omega} (\nabla\overline{\psi}B\psi - \nabla\overline{\Tilde{\psi}}\Tilde{B}\Tilde{\psi})\cdot(u-\Tilde{u})}_{\circled{I-5}}
    \end{aligned}
\end{equation}\medskip

Of course, the next step in this process is to reduce the above equation to a form that is amenable to the Gr\"onwall's inequality, and this is usually achieved by employing H\"older's and Young's inequalities to extract quadratic terms of interest. Since the strong solutions are sufficiently smooth, it is easy to see (see Section 2.5 in \cite{Lions1996MathematicalMechanics} for details) that \circled{I-1} and \circled{I-2} can be handled in a straightforward manner, by extracting terms like $\norm{\rho-\Tilde{\rho}}_{L^2_x}^2$ and $\norm{\sqrt{\rho}(u-\Tilde{u})}_{L^2_x}^2$. In the usual Navier-Stokes setting, these are the only terms. However, here we have additional coupling terms. While \circled{I-4} can be handled in the same way as the previous ones, the other integrals pose problems. For instance, during the analysis of \circled{I-5},  one would encounter a term of the form: \medskip

\begin{equation*}
    \int_{\Omega} B\psi \ (u-\Tilde{u}) \cdot \nabla\overline{\varphi}
\end{equation*}\medskip

\noindent where $\varphi = \psi-\Tilde{\psi}$. Recall that based on our assumptions in Conjecture \ref{weak-strong uniqueness}, we only have $\psi\in L^{\infty}_t H^1_x \cap L^2_t H^2_x$, $B\psi\in L^2_{t,x}$ and $u\in L^{\infty}_t L^2_x \cap L^2_t H^1_x$. It is clear that no splitting of Lebesgue indices (using H\"older's inequality) will be satisfactory: $B\psi$ has to be in $L^2_x$, which means the other two factors can be in something like $L^6_x$ and $L^3_x$, say. Then, the one in $L^6_x$ is extracted as a dissipative factor, and leaves $B\psi$ and the other factor squared. We can then use interpolation to replace the $L^3_x$ factor by $L^2_x$ and $L^6_x$. Once again, we have to extract out another dissipative factor, which results in $\norm{B\psi}_{L^2_x}^4$. Other combinations of Lebesgue indices also lead to these kinds of dead-ends. \medskip

We run into a very similar wall with the corresponding term in \circled{I-3}. Mirroring the above argument, since $B\psi\in L^2_x$, the only option we are left with is to have each of the remaining factors in $L^6_x$. Extracting out a dissipative velocity contribution of the form $\norm{\nabla(u-\Tilde{u})}_{L^2_x}^2$ yields: \medskip

\begin{align*}
    \int_{\Omega} B\psi \ (u+\Tilde{u}) \cdot (u-\Tilde{u}) \overline{\varphi} &\lesssim \norm{B\psi}_{L^2_x} \norm{u+\Tilde{u}}_{L^6_x} \norm{u-\Tilde{u}}_{L^6_x} \norm{\varphi}_{L^6_x} \\
    &\lesssim \norm{B\psi}_{L^2_x}^2 \norm{u+\Tilde{u}}_{H^1_x}^2 \norm{\nabla\varphi}_{L^2_x}^2 + \norm{\nabla(u-\Tilde{u})}_{L^2_x}^2
\end{align*}\medskip

Evidently, this too is inadequate, since when the Gr\"onwall inequality is applied, the time integral has to distribute over $\norm{B\psi}_{L^2_x}^2 \norm{u+\Tilde{u}}_{H^1_x}^2$, which is impossible due to the regularity assumed. \medskip

In conclusion, we note that the standard weak-strong uniqueness proof will not work here, even for strong solutions that are smooth. It appears that Conjecture \ref{weak-strong uniqueness} may not hold true. This naturally leads us to consider slightly more regular weak solutions, as is done in Theorem \ref{local uniqueness}. \bigskip

\section{Weak-moderate uniqueness (Proof of Theorem \ref{local uniqueness})} \label{local uniqueness of weak solutions}

We will use the difference between the weak and moderate solutions as the test functions, so as to obtain a Gr\"onwall-type estimate for the $L^{\infty}_{[0,T]} L^2_x$ norm ($L^{\infty}_{[0,T]} H^1_x$ for the wavefunction) of the difference of two solutions. This is possible in the case of the velocity and wavefunction simply because of their regularities. In the case of the density, it is possible to achieve it by working with smooth approximations and passing to the limit (recall that $\rho\in C^0_t L^2_x$). Finally, observe that the increased spatial regularity of $\Tilde{\rho}$ results in an increased temporal regularity too.

\begin{equation} \label{additional time regularity of density}
\begin{aligned}
    \norm{\partial_t \Tilde{\rho}}_{L^2_{[0,T]} L^2_x} &\lesssim \norm{\Tilde{u}\cdot\nabla\Tilde{\rho}}_{L^2_{[0,T]} L^2_x} + \norm{\overline{\Tilde{\psi}}\Tilde{B}\Tilde{\psi}}_{L^2_{[0,T]} L^2_x} \\
    &\lesssim \norm{\Tilde{u}}_{L^{\infty}_{[0,T]} L^6_x} \norm{\Tilde{\rho}}_{L^2_{[0,T]} W^{1,3}_x} + \norm{\Tilde{\psi}}_{L^{\infty}_{[0,T]} L^{\infty}_x} \norm{\Tilde{B}\Tilde{\psi}}_{L^2_{[0,T]} L^2_x} < \infty
\end{aligned}
\end{equation}\medskip


    
    

The basic idea of the proof is to subtract the governing equations for the weak and moderate fields, and use the difference of the fields as the test functions itself, i.e., $\varphi = \psi - \Tilde{\psi}$, $\Phi = u - \Tilde{u}$, and $\sigma = \rho - \Tilde{\rho}$. In what follows, $\kappa$ is a positive real number used to extract out the dissipative terms $\norm{\nabla\Phi}_{L^2_x}^2$ and $\norm{D^2\varphi}_{L^2_x}^2$, and it is small enough that they may be absorbed into the dissipative terms on the LHS of the momentum and Schr\"odinger equations, respectively. \medskip

\subsection{The wavefunction equation} \label{weak uniqueness for wavefunction}

We will begin by acting with the gradient operator on \eqref{NLS} for each of the two solutions $\psi$ and $\Tilde{\psi}$. Then, we will subtract the equations and integrate against $\nabla(\overline{\psi} - \overline{\Tilde{\psi}})$. The real part of the resulting equation is given below. The analysis is similar to the calculations in Section 3.2 of \cite{Jayanti2022LocalSuperfluidity}, in that we will be looking at the energy ($H^1_x$ norm) difference between the two wavefunctions. \medskip

\begin{equation} \label{wavefunction gradient difference step 1}
    \frac{d}{dt}\frac{1}{2} \norm{\nabla(\psi - \Tilde{\psi})}_{L^2_x}^2 = -\Lambda\Re\int_{\Omega} \nabla(\overline{\psi}-\overline{\Tilde{\psi}})\cdot\nabla(B\psi - \Tilde{B}\Tilde{\psi}) + \mu\Im \int_{\Omega} \nabla(\overline{\psi}-\overline{\Tilde{\psi}})\cdot\nabla\left(\abs{\psi}^2\psi - \abs{\Tilde{\psi}}^2\Tilde{\psi} \right)
\end{equation}\medskip

Integrating by parts on the RHS, we get a Laplacian term in each of the integrals. In the first integral, we observe that $B\psi - \Tilde{B}\Tilde{\psi} = -\frac{1}{2}\Delta(\psi - \Tilde{\psi}) + \dots$, so that this gives us a ``dissipation" term along with other terms. In the second term of the RHS, we extract the $\Delta(\psi - \Tilde{\psi})$ in $L^2_x$, and absorb it into the aforementioned dissipation term. \medskip

\begin{align*}
    -\Lambda\Re\int_{\Omega} \nabla\overline{\varphi}\cdot\nabla(B\psi - \Tilde{B}\Tilde{\psi}) &\begin{multlined}[t]= \Lambda\Re\int_{\Omega} \Delta\varphi \left[-\frac{1}{2}\Delta\varphi + \frac{1}{2}\abs{u}^2\varphi + \frac{1}{2}\Tilde{\psi}(u+\Tilde{u})\cdot\Phi + i\Phi\cdot\nabla\psi \right. \\
    \left. + \Tilde{u}\cdot\nabla\varphi + \mu\left( \abs{\psi}^2\varphi + \Tilde{\psi}\abs{\varphi}^2 + \abs{\Tilde{\psi}}^2\varphi + \Tilde{\psi}^2 \overline{\varphi} \right) \right] \end{multlined} \\
    &\begin{multlined}[t] 
    \lesssim -\frac{\Lambda}{2}\norm{D^2\varphi}_{L^2_x}^2 + \frac{\kappa}{2}\Lambda\norm{D^2\varphi}_{L^2_x}^2 \\ + \kappa^{-1}\Lambda\left\Vert\abs{u}^2\varphi + \Tilde{\psi}(u+\Tilde{u})\cdot\Phi + i\Phi\cdot\nabla\psi + \Tilde{u}\cdot\nabla\varphi \right. \\ \left. + \mu\left( \abs{\psi}^2\varphi + \Tilde{\psi}\abs{\varphi}^2 + \abs{\Tilde{\psi}}^2\varphi + \Tilde{\psi}^2 \overline{\varphi} \right)\right\Vert_{L^2_x}^2 \end{multlined}
\end{align*}
\medskip

The last term in \eqref{wavefunction gradient difference step 1} is also handled in a similar manner, and after absorbing the ``dissipation'' term, we arrive at: \medskip






\begin{multline} \label{varphi equation initial}
    \frac{d}{dt}\norm{\nabla\varphi}_{L^2_x}^2 + \frac{\Lambda}{2}\norm{D^2\varphi}_{L^2_x}^2 \lesssim \kappa^{-1}\Lambda \left[ \underbrace{\norm{\abs{u}^2\varphi}_{L^2_x}^2}_{\circled{A-1}} + \underbrace{\norm{\Tilde{\psi}(u+\Tilde{u})\cdot\Phi}_{L^2_x}^2}_{\circled{A-2}} + \underbrace{\norm{\Phi\cdot\nabla\psi}_{L^2_x}^2}_{\circled{A-3}} + \underbrace{\norm{\Tilde{u}\cdot\nabla\varphi}_{L^2_x}^2}_{\circled{A-4}} \right. \\ \left. + \mu^2 (1+\Lambda^{-1}) \underbrace{\norm{\abs{\psi}^2\varphi + \Tilde{\psi}\abs{\varphi}^2 + \abs{\Tilde{\psi}}^2\varphi + \Tilde{\psi}^2 \overline{\varphi}}_{L^2_x}^2}_{\circled{A-5} + \circled{A-6} + \circled{A-7} + \circled{A-8}} \right] + \kappa\Lambda \norm{D^2\varphi}_{L^2_x}^2
\end{multline}\medskip



Now, we will look at each of the terms in \eqref{varphi equation initial}.

\begin{enumerate}[(i)]
    \item 
    \begin{equation*}
        \circled{A-1} \lesssim \Lambda\norm{u}_{L^6_x}^4\norm{\varphi}_{L^6_x}^2 \lesssim \Lambda\norm{u}_{H^1_x}^4 \norm{\nabla\varphi}_{L^2_x}^2
    \end{equation*}
    
    \item
    \begin{align*}
        \circled{A-2} &\lesssim \kappa^{-1}\Lambda \norm{\Tilde{\psi}}_{L^{\infty}_x}^2 \norm{u+\Tilde{u}}_{L^6_x}^2 \norm{\Phi}_{L^3_x}^2  \\ 
        &\lesssim \kappa^{-3}\nu^{-1}\Lambda^2 \norm{u+\Tilde{u}}_{H^1_x}^4 \norm{\Tilde{\psi}}_{H^2_x}^4 \norm{\Phi}_{L^2_x}^2 + \kappa\nu\norm{\nabla\Phi}_{L^2_x}^2
    \end{align*}
    
    \item\begin{align*}
        \circled{A-3} &\lesssim \kappa^{-1}\Lambda \norm{\nabla\psi}_{L^6_x}^2 \norm{\Phi}_{L^3_x}^2 \\
        &\lesssim \kappa^{-3}\nu^{-1}\Lambda^2 \norm{\psi}_{H^2_x}^4 \norm{\Phi}_{L^2_x}^2 + \kappa\nu\norm{\nabla\Phi}_{L^2_x}^2
    \end{align*}\medskip
    
    \item
    \begin{align*}
        \circled{A-4} &\lesssim \kappa^{-1}\Lambda \norm{\Tilde{u}}_{L^6_x}^2 \norm{\nabla\varphi}_{L^3_x}^2 \\
        &\lesssim \kappa^{-3}\Lambda \norm{\Tilde{u}}_{H^1_x}^4 \norm{\nabla\varphi}_{L^2_x}^2 + \kappa\Lambda\norm{D^2\varphi}_{L^2_x}^2
    \end{align*}\medskip
    
    \noindent where we have used $L^3$ interpolation. Note that the last term above, combined with the last term on the RHS of \eqref{varphi equation initial}, can be absorbed into the LHS of \eqref{varphi equation initial}. This is the ``dissipation'' term we seek to extract out whenever appropriate. \medskip
    
    \item
    \begin{align*}
        \circled{A-5} &\lesssim \kappa^{-1}(1+\Lambda)\mu^2 \norm{\psi}_{L^6_x}^4 \norm{\varphi}_{L^6_x}^2 \\
        &\lesssim \kappa^{-1}(1+\Lambda)\mu^2 \norm{\psi}_{H^1_x}^4 \norm{\nabla\varphi}_{L^2_x}^2
    \end{align*}\medskip
    
    \item
    \begin{align*}
        \circled{A-6} &\lesssim \kappa^{-1}(1+\Lambda)\mu^2 \norm{\Tilde{\psi}}_{L^6_x}^2 \norm{\varphi}_{L^6_x}^2 \norm{\varphi}_{L^6_x}^2 \\
        &\lesssim \kappa^{-1}(1+\Lambda)\mu^2 \norm{\Tilde{\psi}}_{H^1_x}^2 \left( \norm{\psi}_{H^1_x}^2 + \norm{\Tilde{\psi}}_{H^1_x}^2 \right) \norm{\nabla\varphi}_{L^2_x}^2
    \end{align*}\medskip
    
    \item
    \begin{align*}
        \circled{A-7} + \circled{A-8} &\lesssim \kappa^{-1}(1+\Lambda)\mu^2 \norm{\Tilde{\psi}}_{L^6_x}^4 \norm{\varphi}_{L^6_x}^2 \\
        &\lesssim \kappa^{-1}(1+\Lambda)\mu^2 \norm{\Tilde{\psi}}_{H^1_x}^4 \norm{\nabla\varphi}_{L^2_x}^2
    \end{align*}\medskip
    
\end{enumerate}

Combining these estimates into \eqref{varphi equation initial} gives us:

\begin{equation} \label{varphi equation final}
    \frac{d}{dt} \norm{\nabla\varphi}_{L^2_x}^2 + \frac{\Lambda}{2} \norm{D^2\varphi}_{L^2_x}^2 \le h_1(t) \left[ \norm{\nabla\varphi}_{L^2_x}^2 + \norm{\Phi}_{L^2_x}^2 \right] + \kappa\nu\norm{\nabla\Phi}_{L^2_x}^2 + \kappa\Lambda\norm{D^2\varphi}_{L^2_x}^2
\end{equation}\medskip

\noindent where $h_1 \in L^1_{[0,T]}$. \bigskip

\subsection{The momentum equation} \label{the momentum equation uniqueness}

The notion of the weak solution holds for all times $0\le t\le T$, where $T$ is the chosen local existence time. Considering  \eqref{weak solution velocity} and \eqref{weak solution density} up to a time $t \le T$, and differentiating with respect to $t$, we arrive at:

    
\begin{equation} \label{weak solution velocity differentiated wrt time}
    -\int_{\Omega} \left[ \rho u\cdot \partial_t \Phi + \rho u\otimes u:\nabla\Phi - \nu\nabla u:\nabla\Phi - 2\Lambda\Phi\cdot\Im(\nabla\Bar{\psi}B\psi) \right] = -\frac{d}{dt}\int_{\Omega} \rho(t)u(t)\cdot\Phi(t)
\end{equation}\medskip
    
\begin{equation} \label{weak solution density differentiated wrt time}
    -\int_{\Omega} \left[ \rho \partial_t \sigma + \rho u\cdot\nabla\sigma + 2\Lambda \sigma \Re(\Bar{\psi}B\psi) \right] = -\frac{d}{dt}\int_{\Omega} \rho(t)\sigma(t) 
\end{equation}\medskip

Subtracting \eqref{weak solution velocity differentiated wrt time} written for $u$ and $\Tilde{u}$, we get:

\begin{equation} \label{weak solution velocity difference}
\begin{multlined} 
    -\int_{\Omega} \rho \Phi\cdot \partial_t \Phi - \int_{\Omega}\sigma \Tilde{u}\cdot\partial_t \Phi - \int_{\Omega} \left[ \rho\Phi\otimes u + \rho \Tilde{u}\otimes\Phi + \sigma\Tilde{u}\otimes\Tilde{u} \right]:\nabla\Phi + \nu\norm{\nabla\Phi}_{L^2_x}^2 \\ + 2\Lambda\int_{\Omega} \Phi\cdot \Im\left( \nabla\overline{\psi}B\psi - \nabla\overline{\Tilde{\psi}}\Tilde{B}\Tilde{\psi} \right) = -\frac{d}{dt}\int_{\Omega} \rho\abs{\Phi}^2 - \frac{d}{dt} \int_{\Omega} \sigma \Tilde{u}\cdot\Phi
\end{multlined}
\end{equation}\medskip

The first term on the LHS can be rewritten since we see that $\frac{\abs{\Phi}^2}{2}$ satisfies the requirements to be a test function for the density field in \eqref{weak solution density differentiated wrt time}.

\begin{equation*}
    -\int_{\Omega} \rho \Phi\cdot \partial_t \Phi = -\int_{\Omega} \rho \partial_t \frac{\abs{\Phi}^2}{2} = \int_{\Omega} \rho u\cdot \nabla \frac{\abs{\Phi}^2}{2} + \Lambda\Re\int_{\Omega} \abs{\Phi}^2 \overline{\psi}B\psi - \frac{d}{dt}\int_{\Omega} \rho \frac{\abs{\Phi}^2}{2}
\end{equation*}\medskip

Further simplification using H\"older's and Young's inequalities (to appropriately absorb the dissipation $\norm{\nabla\Phi}_{L^2_x}^2$ from the RHS into the LHS) yields,

\begin{equation} \label{phi equation initial}
    \begin{aligned}
        \renewcommand{\arraystretch}{2}
        \left. \begin{array}{c}
         \frac{d}{dt}\frac{1}{2}\norm{\sqrt{\rho}\Phi}_{L^2_x}^2  \\
          + \ \frac{\nu}{2}\norm{\nabla\Phi}_{L^2_x}^2 
        \end{array}  \right\} \le &\frac{1}{\nu} \underbrace{\left[ \norm{\rho\Phi\otimes u}_{L^2_x}^2 + \norm{\rho \Tilde{u}\otimes\Phi}_{L^2_x}^2 + \norm{\sigma\Tilde{u}\otimes\Tilde{u}}_{L^2_x}^2 \right]}_{\circled{B-1} + \circled{B-2} + \circled{B-3}} \\
        &- 2\Lambda \underbrace{\int_{\Omega} \Phi\cdot\Im\left( \nabla\overline{\varphi}B\psi + \nabla\overline{\Tilde{\psi}}B\varphi + \nabla\overline{\Tilde{\psi}}(B-\Tilde{B})\Tilde{\psi} \right) }_{\circled{B-4} + \circled{B-5} + \circled{B-6}} \\
        &- \underbrace{\int_{\Omega}\partial_t (\sigma \Tilde{u})\cdot\Phi}_{\circled{B-7}} - \underbrace{\int_{\Omega}\rho u\cdot \nabla\frac{\abs{\Phi}^2}{2}}_{\circled{B-8}} - \Lambda\underbrace{\int_{\Omega}\abs{\Phi}^2 \Re\overline{\psi}B\psi }_{\circled{B-9}}
    \end{aligned}
\end{equation}\medskip

Each of the terms on the RHS can be bounded from above as follows (using the H\"older and Young inequalities, and the results listed in Lemmas \ref{sobolev embeddings} and \ref{inequalities and interpolations}).

\begin{enumerate}[(i)]
    \item 
    \begin{equation*}
        \circled{B-1} \lesssim \nu^{-1}\norm{\rho}_{L^{\infty}_x}^2 \norm{u}_{L^{\infty}_x}^2 \norm{\Phi}_{L^2_x}^2 \lesssim \nu^{-1}\norm{\rho}_{L^{\infty}_x}^2 \norm{u}_{H^2_x}^2 \norm{\Phi}_{L^2_x}^2
    \end{equation*}\medskip
    
    \item 
    \begin{equation*}
        \circled{B-2} \lesssim \nu^{-1}\norm{\rho}_{L^{\infty}_x}^2 \norm{\Tilde{u}}_{L^{\infty}_x}^2 \norm{\Phi}_{L^2_x}^2 \lesssim \nu^{-1}\norm{\rho}_{L^{\infty}_x}^2 \norm{\Tilde{u}}_{H^2_x}^2 \norm{\Phi}_{L^2_x}^2
    \end{equation*}\medskip
    
    \item Using Agmon's inequality to bound the $L^{\infty}$ norm,
    \begin{equation*}
        \circled{B-3} \lesssim \nu^{-1}\norm{\Tilde{u}}_{L^{\infty}_x}^4 \norm{\sigma}_{L^2_x}^2 \lesssim \nu^{-1}\norm{\Tilde{u}}_{H^1_x}^2 \norm{\Tilde{u}}_{H^2_x}^2 \norm{\sigma}_{L^2_x}^2
    \end{equation*}\medskip
    
    \item We interpolate the $L^3$ norm (Lemma \ref{inequalities and interpolations}):
    \begin{align*}
        \circled{B-4} &\lesssim \Lambda\norm{\nabla\varphi}_{L^6_x} \norm{B\psi}_{L^3_x} \norm{\Phi}_{L^2_x} \\
        &\lesssim \kappa^{-1}\Lambda \norm{B\psi}_{L^3_x}^2 \norm{\Phi}_{L^2_x}^2 + \kappa\Lambda \norm{D^2\varphi}_{L^2_x}^2 \\
        &\lesssim \kappa^{-1}\Lambda \norm{B\psi}_{L^2_x}\norm{B\psi}_{H^1_x} \norm{\Phi}_{L^2_x}^2 + \kappa\Lambda \norm{D^2\varphi}_{L^2_x}^2
    \end{align*}\medskip
    
    


    \item
    \begin{align*}
        \circled{B-5} &= \Lambda\Im\int_{\Omega} \Phi\cdot\nabla\overline{\Tilde{\psi}} \left( -\frac{1}{2}\Delta\varphi + \frac{1}{2}\abs{u}^2\varphi + iu\cdot\nabla\varphi + \mu\abs{\psi}^2\varphi \right) \\ 
        &= \circled{B-5.1} + \circled{B-5.2} + \circled{B-5.3} + \circled{B-5.4}
    \end{align*}\medskip
    
    \begin{itemize}
        \item
    \begin{align*}
        \circled{B-5.1} &\lesssim \Lambda\norm{\Phi}_{L^3_x} \norm{\nabla\Tilde{\psi}}_{L^6_x} \norm{\Delta\varphi}_{L^2_x} \\
        &\lesssim \kappa^{-1}\Lambda \norm{\Tilde{\psi}}_{H^2_x}^2 \norm{\Phi}_{L^3_x}^2 + \kappa\Lambda\norm{\Delta\varphi}_{L^2_x}^2 \\
        &\lesssim \kappa^{-3}\nu^{-1}\Lambda \norm{\Tilde{\psi}}_{H^2_x}^4 \norm{\Phi}_{L^2_x}^2 + \kappa\nu\norm{\nabla\Phi}_{L^2_x}^2 + \kappa\Lambda\norm{D^2\varphi}_{L^2_x}^2
    \end{align*}\medskip
    
    

    
    \item 
    \begin{align*}
        \circled{B-5.2} &\lesssim \Lambda\norm{u}_{L^6_x}^2 \norm{\nabla\Tilde{\psi}}_{L^6_x} \norm{\Phi}_{L^3_x} \norm{\varphi}_{L^6_x} \\ 
        &\lesssim (\kappa\nu)^{-1}\norm{u}_{H^1_x}^4 \norm{\Phi}_{L^2_x}^2 + \Lambda^2 \norm{\Tilde{\psi}}_{H^2_x}^2 \norm{\varphi}_{L^2_x}^2 + \kappa\nu \norm{\nabla\Phi}_{L^2_x}^2
    \end{align*}\medskip
    
    \item
    \begin{align*}
        \circled{B-5.3} &\lesssim \Lambda\norm{u}_{L^6_x} \norm{\nabla\Tilde{\psi}}_{L^6_x} \norm{\Phi}_{L^6_x} \norm{\nabla\varphi}_{L^2_x} \\ 
        &\lesssim (\kappa\nu)^{-1}\Lambda^2\norm{u}_{H^1_x}^2 \norm{\Tilde{\psi}}_{H^2_x}^2 \norm{\nabla\varphi}_{L^2_x}^2 + \kappa\nu\norm{\nabla\Phi}_{L^2_x}^2
    \end{align*}\medskip
    
    \item
    \begin{align*}
        \circled{B-5.4} &\lesssim \Lambda \norm{\psi}_{L^{\infty}_x}^2\norm{\nabla\Tilde{\psi}}_{L^6_x} \norm{\Phi}_{L^3_x} \norm{\varphi}_{L^6_x} \\
        &\lesssim (\kappa\nu)^{-1} \norm{\Tilde{\psi}}_{H^2_x}^4 \norm{\Phi}_{L^2_x}^2 + \Lambda^2\norm{\psi}_{H^2_x}^4 \norm{\nabla\varphi}_{L^2_x}^2 + \kappa\nu\norm{\nabla\Phi}_{L^2_x}^2
    \end{align*}\medskip
    
    \end{itemize}

    \item First we expand $B-\Tilde{B}$ as: \medskip
    
    \begin{align} \label{B - Tilde(B)}
        B - \Tilde{B} &= \frac{1}{2}\left( \abs{u}^2 - \abs{\Tilde{u}}^2 \right) + i(u-\Tilde{u})\cdot\nabla + \mu\left( \abs{\psi}^2 - \abs{\Tilde{\psi}}^2 \right) \notag \\
        &= \frac{1}{2} (u + \Tilde{u})\cdot\Phi + i\Phi\cdot\nabla + \mu \left( \abs{\varphi}^2 + \overline{\varphi}\Tilde{\psi} + \varphi\overline{\Tilde{\psi}} \right)
    \end{align}\medskip

    Then,

    \begin{align*}
        \circled{B-6} &= \int_{\Omega} \Phi\cdot\Im \left[ \nabla\overline{\Tilde{\psi}} \left\{ \frac{1}{2}(u+\Tilde{u})\cdot\Phi + i\Phi\cdot\nabla + \mu(\abs{\varphi}^2 + \overline{\varphi}\Tilde{\psi} + \varphi \overline{\Tilde{\psi}})\Tilde{\psi} \right\} \Tilde{\psi} \right] \\
        &\begin{multlined} \lesssim \norm{\Tilde{\psi}}_{L^{\infty}_x} \norm{\nabla\Tilde{\psi}}_{L^6_x} \norm{u+\Tilde{u}}_{L^6_x} \norm{\Phi}_{L^3_x}^2 + \norm{\nabla\Tilde{\psi}}_{L^6_x}^2 \norm{\Phi}_{L^3_x}^2 + \norm{\Tilde{\psi}}_{L^{\infty}_x}^3 \norm{\nabla\Tilde{\psi}}_{L^3_x} \norm{\Phi}_{L^2_x} \norm{\varphi}_{L^6_x} \\ + \norm{\Tilde{\psi}}_{L^{\infty}_x}^2 \norm{\nabla\Tilde{\psi}}_{L^6_x} \norm{\Phi}_{L^2_x} \left(\norm{\psi}_{L^6_x} + \norm{\Tilde{\psi}}_{L^6_x}\right) \norm{\varphi}_{L^6_x} \end{multlined} \\
        &\begin{multlined} \lesssim (\kappa\nu)^{-1}\norm{\Tilde{\psi}}_{H^2_x}^4 \left( 1 + \norm{u+\Tilde{u}}_{H^1_x}^2 \right) \norm{\Phi}_{L^2_x}^2 + \norm{\Tilde{\psi}}_{H^2_x}^2 \left( \norm{\nabla\psi}_{L^2_x}^2 + \norm{\nabla\Tilde{\psi}}_{L^2_x}^2 \right) \norm{\nabla\varphi}_{L^2_x}^2 \\ + \kappa\nu\norm{\nabla\Phi}_{L^2_x}^2 \end{multlined}
    \end{align*}\medskip
    
    \item 
    \begin{equation*}
        \circled{B-7} = -\underbrace{\int_{\Omega}\partial_t \Tilde{u} \cdot\sigma\Phi}_{\circled{B-7.1}} - \underbrace{\int_{\Omega} \partial_t \sigma (\Tilde{u}\cdot\Phi)}_{\circled{B-7.2}}
    \end{equation*}\medskip

    \begin{itemize}
    \item
    
    \begin{align*}
        \circled{B-7.1} &\lesssim \norm{\partial_t \Tilde{u}}_{L^3_x} \norm{\sigma}_{L^2_x} \norm{\Phi}_{L^6_x} \\
        &\lesssim (\kappa\nu)^{-1} \norm{\partial_t \Tilde{u}}_{L^3_x}^2 \norm{\sigma}_{L^2_x}^2 + \kappa\nu\norm{\nabla\Phi}_{L^2_x}^2 
    \end{align*}\medskip
    
    \item By subtracting the continuity equations for $\rho$ and $\Tilde{\rho}$, we obtain 
    
    \begin{equation*}
        \partial_t \sigma + \Phi\cdot\nabla\Tilde{\rho} +u\cdot\nabla\sigma = 2\Lambda \Re \left( \overline{\psi}B\psi - \overline{\Tilde{\psi}}\Tilde{B}\Tilde{\psi} \right)
    \end{equation*}
    
    which allows us to rewrite $\circled{B-7.2}$ as
    
    \begin{align*}
        \circled{B-7.2} = \underbrace{\int_{\Omega} \left(\Tilde{u}\cdot\Phi\right) \left[ -\Phi\cdot\nabla\Tilde{\rho} - u\cdot\nabla\sigma + 2\Lambda \Re \left( \overline{\psi}B\psi - \overline{\Tilde{\psi}}\Tilde{B}\Tilde{\psi} \right) \right]}_{\circled{B-7.2.1}+\circled{B-7.2.2}+\circled{B-7.2.3}}
    \end{align*}\medskip
    
    Each of these terms will now be estimated. For the first two terms, we will integrate by parts and interpolate the $L^3$ norm: \medskip
    
    \begin{align*}
        \circled{B-7.2.1} &= \int_{\Omega} \Tilde{\rho}\Phi\cdot\nabla(\Phi\cdot\Tilde{u}) = \int_{\Omega} \Tilde{\rho}\Phi\cdot(\nabla\Phi)\cdot\Tilde{u} + \int_{\Omega} \Tilde{\rho}\Phi\otimes\Phi:\nabla\Tilde{u}) \\
        &\lesssim \norm{\Tilde{\rho}}_{L^{\infty}_x} \norm{\Phi}_{L^3_x} \norm{\nabla\Phi}_{L^2_x} \norm{\Tilde{u}}_{L^6_x} + \norm{\Tilde{\rho}}_{L^{\infty}_x} \norm{\Phi}_{L^3_x} \norm{\Phi}_{L^6_x} \norm{\nabla\Tilde{u}}_{L^2_x} \\
        &\lesssim (\kappa\nu)^{-1}\norm{\Tilde{\rho}}_{L^{\infty}_x}^2 \norm{\Tilde{u}}_{H^1_x}^2 \norm{\Phi}_{L^3_x}^2 + \kappa\nu\norm{\nabla\Phi}_{L^2_x}^2 \\
        &\lesssim (\kappa\nu)^{-3}\norm{\Tilde{\rho}}_{L^{\infty}_x}^4 \norm{\Tilde{u}}_{H^1_x}^4 \norm{\Phi}_{L^2_x}^2 + \kappa\nu\norm{\nabla\Phi}_{L^2_x}^2
    \end{align*} \medskip
    
    
    \begin{align*}
        \circled{B-7.2.2} &= \int_{\Omega}\sigma u\cdot\nabla(\Phi\cdot\Tilde{u}) = \int_{\Omega}\sigma u\cdot(\nabla\Phi)\cdot\Tilde{u} + \int_{\Omega}\sigma u\cdot(\nabla\Tilde{u})\cdot\Phi \\
        &\lesssim \norm{u}_{L^{\infty}_x} \norm{\Tilde{u}}_{L^{\infty}_x} \norm{\sigma}_{L^2_x} \norm{\nabla\Phi}_{L^2_x} + \norm{u}_{L^6_x} \norm{\nabla\Tilde{u}}_{L^6_x} \norm{\sigma}_{L^2_x} \norm{\Phi}_{L^6_x} \\
        &\lesssim (\kappa\nu)^{-1}\left( \norm{u}_{L^{\infty}_x}^2 \norm{\Tilde{u}}_{L^{\infty}_x}^2 + \norm{u}_{H^1_x}^2 \norm{\Tilde{u}}_{H^2_x}^2 \right) \norm{\sigma}_{L^2_x}^2 + \kappa\nu\norm{\nabla\Phi}_{L^2_x}^2 \\
        &\begin{multlined} \lesssim (\kappa\nu)^{-1}\left( \norm{u}_{H^1_x} \norm{u}_{H^2_x} \norm{\Tilde{u}}_{H^1_x} \norm{\Tilde{u}}_{H^2_x} + \norm{u}_{H^1_x}^2 \norm{\Tilde{u}}_{H^2_x}^2 \right) \norm{\sigma}_{L^2_x}^2 + \kappa\nu\norm{\nabla\Phi}_{L^2_x}^2 \end{multlined}
    \end{align*} \medskip
    
    \noindent where we used Agmon's inequality to bound the $L^{\infty}$ terms. Finally,
    
    \begin{align} \label{B-7.2.3}
        \circled{B-7.2.3} = \underbrace{2\Lambda\int_{\Omega} (\Tilde{u}\cdot\Phi) \Re \left[ \overline{\varphi}\Tilde{B}\Tilde{\psi} + \overline{\psi}\Tilde{B}\varphi + \overline{\psi} (B-\Tilde{B}) \psi \right]}_{\circled{B-7.2.3.1}+\circled{B-7.2.3.2}+\circled{B-7.2.3.3}}
    \end{align}\medskip
    
    These are, in turn, estimated as follows. \medskip
    
    \begin{align*}
        \circled{B-7.2.3.1} &\lesssim \Lambda\norm{\Tilde{u}}_{L^6_x} \norm{\Phi}_{L^6_x} \norm{\varphi}_{L^6_x} \norm{\Tilde{B}\Tilde{\psi}}_{L^2_x} \\
        &\lesssim (\kappa\nu)^{-1}\Lambda^2\norm{\Tilde{u}}_{H^1_x}^2 \norm{\Tilde{B}\Tilde{\psi}}_{L^2_x}^2 \norm{\nabla\varphi}_{L^2_x}^2 + \kappa\nu\norm{\nabla\Phi}_{L^2_x}^2
    \end{align*}\medskip
    
    The middle term requires to be broken down further, much like $\circled{B-5}$.
    
    \begin{align*}
        \circled{B-7.2.3.2} &= \underbrace{\Lambda\Re\int_{\Omega} (\Tilde{u}\cdot\Phi) \overline{\psi} \left[ -\frac{1}{2}\Delta\varphi + \frac{1}{2}\abs{\Tilde{u}}^2\varphi + i\Tilde{u}\cdot\nabla\varphi + \mu \abs{\Tilde{\psi}}^2 \varphi \right]}_{\circled{B-7.2.3.2.1}+\circled{B-7.2.3.2.2}+\circled{B-7.2.3.2.3}+\circled{B-7.2.3.2.4}}
    \end{align*}\medskip
    
    Just as in \circled{B-5}, we extract a dissipative factor, and make use of interpolation. \medskip
    
    \begin{align*}
        \circled{B-7.2.3.2.1} &\lesssim \Lambda\norm{\Tilde{u}}_{L^6_x} \norm{\Phi}_{L^3_x} \norm{\psi}_{L^{\infty}_x} \norm{\Delta\varphi}_{L^2_x} \\
        &\lesssim \kappa^{-1}\Lambda \norm{\Tilde{u}}_{H^1_x}^2 \norm{\psi}_{H^2_x}^2 \norm{\Phi}_{L^2_x} \norm{\Phi}_{L^6_x} + \kappa\Lambda \norm{D^2\varphi}_{L^2_x}^2 \\
        &\lesssim \kappa^{-3}\nu^{-1}\Lambda^2 \norm{\Tilde{u}}_{H^1_x}^4 \norm{\psi}_{H^2_x}^4 \norm{\Phi}_{L^2_x}^2 + \kappa\nu\norm{\nabla\Phi}_{L^2_x}^2 + \kappa\Lambda \norm{D^2\varphi}_{L^2_x}^2
    \end{align*}
    
    
    \begin{align*}
        \circled{B-7.2.3.2.2} &\lesssim \Lambda\norm{\Tilde{u}}_{L^6_x}^3 \norm{\psi}_{L^{\infty}_x} \norm{\varphi}_{L^6_x} \norm{\Phi}_{L^6_x} \\
        &\lesssim (\kappa\nu)^{-1}\Lambda^2 \norm{\Tilde{u}}_{H^1_x}^6 \norm{\psi}_{H^2_x}^2 \norm{\nabla\varphi}_{L^2_x}^2 + \kappa\nu\norm{\nabla\Phi}_{L^2_x}^2
    \end{align*}
    
    \begin{align*}
        \circled{B-7.2.3.2.3} &\lesssim \Lambda\norm{\Tilde{u}}_{L^6_x}^2 \norm{\psi}_{L^{\infty}_x} \norm{\nabla\varphi}_{L^6_x} \norm{\Phi}_{L^2_x} \\
        &\lesssim \kappa^{-1}\Lambda \norm{\Tilde{u}}_{H^1_x}^4 \norm{\psi}_{H^2_x}^2 \norm{\Phi}_{L^2_x}^2 + \kappa\Lambda\norm{D^2\varphi}_{L^2_x}^2
    \end{align*}
    
    \begin{align*}
        \circled{B-7.2.3.2.4} &\lesssim \Lambda \norm{\Tilde{u}}_{L^6_x} \norm{\psi}_{L^6_x} \norm{\Tilde{\psi}}_{L^{\infty}_x}^2 \norm{\Phi}_{L^2_x} \norm{\varphi}_{L^6_x} \\
        &\lesssim \Lambda \norm{\Tilde{u}}_{H^1_x}^2 \norm{\psi}_{H^1_x}^2 \norm{\Phi}_{L^2_x}^2 + \Lambda \norm{\Tilde{\psi}}_{H^2_x}^4 \norm{\nabla\varphi}_{L^2_x}^2
    \end{align*}\medskip
    
    The final term in \eqref{B-7.2.3} is dealt with in a manner mirroring that of $\circled{B-6}$.
    
    \begin{align*}
        \circled{B-7.2.3.3} &= 2\Lambda\Re\int_{\Omega} \Tilde{u}\cdot\Phi \overline{\psi} \left[ \frac{1}{2}\Phi\cdot(u+\Tilde{u}) \psi + i\Phi\cdot\nabla\psi + \mu\left( 2\Re(\Tilde{\psi}\overline{\varphi}) + \abs{\varphi}^2 \right)\psi \right] \\
        &\begin{multlined} \lesssim \Lambda \norm{\Tilde{u}}_{L^6_x} \left[ \norm{u+\Tilde{u}}_{L^6_x} \norm{\psi}_{L^{\infty}_x}^2 + \norm{\psi}_{L^{\infty}_x} \norm{\nabla\psi}_{L^6_x} \right] \norm{\Phi}_{L^3_x}^2 \\ + \Lambda \norm{\Tilde{u}}_{L^6_x} \norm{\psi}_{L^{\infty}_x}^2 \left[ \mu \norm{\Tilde{\psi}}_{L^6_x} + \mu\norm{\varphi}_{L^6_x} \right] \norm{\varphi}_{L^6_x} \norm{\Phi}_{L^2_x} \end{multlined} \\
        &\begin{multlined} \lesssim (\kappa\nu)^{-1}\Lambda^2 \norm{\Tilde{u}}_{H^1_x}^2 \norm{\psi}_{H^2_x}^4 \left[ 1+\norm{u+\Tilde{u}}_{H^1_x}^2 + \kappa\nu \right] \norm{\Phi}_{L^2_x}^2 \\ + \mu^2\left[ \norm{\Tilde{\psi}}_{H^1_x}^2 + \norm{\psi}_{H^1_x}^2 \right] \norm{\nabla\varphi}_{L^2_x}^2 + \kappa\nu\norm{\nabla\Phi}_{L^2_x}^2
        \end{multlined}
    \end{align*}\medskip
    
    \end{itemize}

    \item Using \eqref{additional time regularity of density}, and the vector identity $\nabla\frac{\abs{\Phi}^2}{2} = \Phi\cdot\nabla\Phi - (\nabla\times\Phi)\times\Phi$, 
    
    \begin{align*}
        \circled{B-8} &\lesssim \norm{\rho}_{L^{\infty}_x} \norm{u}_{L^6_x} \norm{\Phi}_{L^3_x} \norm{\nabla\Phi}_{L^2_x} \\
        &\lesssim (\kappa\nu)^{-1}\norm{\rho}_{L^{\infty}_x}^2 \norm{u}_{H^1_x}^2 \norm{\Phi}_{L^3_x}^2 + \kappa\nu\norm{\nabla\Phi}_{L^2_x}^2 \\
        &\lesssim (\kappa\nu)^{-3}\norm{\rho}_{L^{\infty}_x}^4 \norm{u}_{H^1_x}^4 \norm{\Phi}_{L^2_x}^2 + \kappa\nu\norm{\nabla\Phi}_{L^2_x}^2
    \end{align*}\medskip
    
    \item 
    \begin{align*}
        \circled{B-9} &\lesssim \norm{\psi}_{L^{\infty}_x} \norm{B\psi}_{L^3_x} \norm{\Phi}_{L^2_x} \norm{\Phi}_{L^6_x} \\
        &\lesssim (\kappa\nu)^{-1}\norm{\psi}_{H^2_x}^2 \norm{B\psi}_{L^2_x} \norm{B\psi}_{H^1_x} \norm{\Phi}_{L^2_x}^2 + \kappa\nu\norm{\nabla\Phi}_{L^2_x}^2
    \end{align*}\medskip
    
\end{enumerate}

Returning to \eqref{phi equation initial}, and using the ensuing estimates, we arrive at:

\begin{equation} \label{phi equation final}
    \varepsilon\frac{d}{dt}\norm{\Phi}_{L^2_x}^2 + \frac{\nu}{2} \norm{\nabla\Phi}_{L^2_x}^2 \le h_2(t)  \left[ \norm{\nabla\varphi}_{L^2_x}^2 + \norm{\Phi}_{L^2_x}^2 + \norm{\sigma}_{L^2_x}^2 \right] + \kappa\nu \norm{\nabla\Phi}_{L^2_x}^2 + \kappa\Lambda \norm{D^2\varphi}_{L^2_x}^2
\end{equation}\medskip

\noindent where $h_2 \in L^1_{[0,T]}$. Also, note that we have replaced the density in the first term of the LHS by its minimum value. \medskip

\subsection{The continuity equation} \label{weak uniqueness for density}

Finally, in the case of the continuity equation, we take the difference of \eqref{continuity} written for each of the solutions, multiply by $\rho-\Tilde{\rho}$ and integrate over $\Omega$. \medskip 



\begin{equation} \label{sigma equation initial}
    \frac{d}{dt} \frac{1}{2} \norm{\sigma}_{L^2_x}^2 = \underbrace{\int_{\Omega} \Tilde{\rho}\Phi\cdot\nabla\sigma}_{\circled{C-1}} - 2\Lambda\Re\underbrace{\int_{\Omega} \sigma\overline{\varphi}B\psi + \sigma\overline{\Tilde{\psi}}(B-\Tilde{B})\Tilde{\psi}}_{\circled{C-2} + \circled{C-3}}
\end{equation}\medskip

\begin{rem}
    The above calculations can be rigorously justified by considering the difference between the equations for $\rho^N$ (the approximate densities in Section 4.3 of \cite{Jayanti2022LocalSuperfluidity}) and $\Tilde{\rho}$. Passing to the limit $N\rightarrow\infty$ leaves us with the desired equation for $\sigma = \rho-\Tilde{\rho}$, since we know that $\rho\in C(0,T;L^2_x)$.
\end{rem}\medskip

\begin{enumerate}[(i)]
    \item Once again, recalling that we have assumed in Theorem \ref{local uniqueness} that $\Tilde{\rho}\in L^2_{[0,T]}W^{1,3}_x$. Thus,
    
    \begin{equation*}
        \circled{C-1} = -\int_{\Omega}\nabla\Tilde{\rho}\cdot\Phi\sigma \lesssim \norm{\nabla\Tilde{\rho}}_{L^3_x} \norm{\Phi}_{L^6_x} \norm{\sigma}_{L^2_x} \lesssim (\kappa\nu)^{-1}\norm{\nabla\Tilde{\rho}}_{L^3_x}^2 \norm{\sigma}_{L^2_x}^2 + \kappa\nu\norm{\nabla\Phi}_{L^2_x}^2
    \end{equation*}\medskip
    
    \noindent Observe that $\sigma$ should always be considered in the $L^2_x$ norm while using H\"older's, since there is no dissipation term for the density. \medskip
    
    \item For the second term, apart from the above, we also interpolate the $L^3_x$ norm between the $L^2_x$ and $L^6_x$ norms.
    
    \begin{align*}
        \circled{C-2} \lesssim \Lambda\norm{\sigma}_{L^2_x} \norm{\varphi}_{L^6_x} \norm{B\psi}_{L^3_x} \lesssim \Lambda^2 \norm{B\psi}_{L^2_x} \norm{\nabla B\psi}_{L^2_x} \norm{\sigma}_{L^2_x}^2 + \norm{\nabla\varphi}_{L^2_x}^2
    \end{align*}\medskip
    
    
    \item Similar analysis on the third term yields
    
    \begin{align*}
        \circled{C-3} &\lesssim \Lambda\norm{\sigma}_{L^2_x} \norm{\Tilde{\psi}}_{L^{\infty}_x} \norm{(B-\Tilde{B})\Tilde{\psi}}_{L^2_x} \\
        &\begin{multlined} \lesssim \Lambda\norm{\sigma}_{L^2_x} \norm{\Tilde{\psi}}_{H^2_x} \left[ \norm{\Tilde{\psi}}_{L^{\infty}_x} \norm{u+\Tilde{u}}_{L^6_x} \norm{\Phi}_{L^3_x} + \norm{\nabla\Tilde{\psi}}_{L^6_x} \norm{\Phi}_{L^3_x} \right. \\ \left. + \norm{\Tilde{\psi}}_{L^6_x}^2 \norm{\varphi}_{L^6_x} + \norm{\Tilde{\psi}}_{L^{\infty}_x} \norm{\varphi}_{L^3_x} \norm{\varphi}_{L^6_x} \right] \end{multlined} \\
        &\begin{multlined} \lesssim \Lambda\left[ \norm{\Tilde{\psi}}_{H^2_x}^4 + \norm{\Tilde{\psi}}_{H^1_x}^2\norm{\Tilde{\psi}}_{H^2_x}^2 \right] \norm{\sigma}_{L^2_x}^2 + (\kappa\nu)^{-1}\Lambda\left[ 1 + \norm{u+\Tilde{u}}_{H^1_x}^4 \right]\norm{\Phi}_{L^2_x}^2 \\
        + \Lambda\left[ \norm{\Tilde{\psi}}_{H^1_x}^2 + \norm{\psi}_{H^1_x}^2 \right] \norm{\nabla\varphi}_{L^2_x}^2 + \kappa\nu \norm{\nabla\Phi}_{L^2_x}^2 \end{multlined}
    \end{align*}\medskip
    
    In going to the last step, we tackle the $\norm{\varphi}_{L^3_x}^2$ term using interpolation and the Poincar\'e inequality as follows:
    
    \begin{equation*}
        \norm{\varphi}_{L^3_x}^2 \lesssim \norm{\varphi}_{L^2_x}\norm{\varphi}_{L^6_x} \lesssim \norm{\varphi}_{L^2_x} \norm{\varphi}_{H^1_x}
        \lesssim \norm{\nabla\varphi}_{L^2_x}^2 \lesssim \norm{\nabla\psi}_{L^2_x}^2 + \norm{\nabla\Tilde{\psi}}_{L^2_x}^2
    \end{equation*}\bigskip

\end{enumerate}

Using the above estimates, we can simplify \eqref{sigma equation initial} to read: \medskip

\begin{equation} \label{sigma equation final}
    \frac{d}{dt}\norm{\sigma}_{L^2_x}^2 \le h_3(t) \norm{\sigma}_{L^2_x}^2 + \kappa\nu\norm{\nabla\Phi}_{L^2_x}^2
\end{equation}\medskip

\noindent where $h_3 \in L^1_{[0,T]}$. \bigskip

\bigskip

\bigskip

We now proceed to add \eqref{varphi equation final}, \eqref{phi equation final} and \eqref{sigma equation final}. By choosing $\kappa$ sufficiently small, we can ensure all the dissipation terms are absorbed by the LHS. What remains is:

\begin{equation}
    \frac{d}{dt} \left[ \norm{\varphi}_{L^2_x}^2 + \norm{\Phi}_{L^2_x}^2 + \norm{\sigma}_{L^2_x}^2 \right] \lesssim \left( h_1 + h_2 + h_3 \right)(t) \left[ \norm{\varphi}_{L^2_x}^2 + \norm{\Phi}_{L^2_x}^2 + \norm{\sigma}_{L^2_x}^2 \right]
\end{equation}\medskip

Using Gr\"onwall's inequality, and the fact that $\norm{\varphi}_{L^2_x} = \norm{\Phi}_{L^2_x} = \norm{\sigma}_{L^2_x} = 0$ at $t=0$ completes the proof of Theorem \ref{local uniqueness}.

\qed

\addtocontents{toc}{\protect\setcounter{tocdepth}{0}}

\section*{Acknowledgments}
Both the authors are grateful to the anonymous referees for their comments and suggestions which significantly helped to improve the original manuscript. P.C.J. was partially supported by the Ann Wylie Fellowship at UMD. Both P.C.J. and K.T. gratefully acknowledge the support of the National Science Foundation under the awards DMS-1614964 and DMS-2008568. \bigskip

\addtocontents{toc}{\protect\setcounter{tocdepth}{2}}



\bibliographystyle{alpha}
\bibliography{references}

%
%

%
%

\end{document}